\documentclass[11pt]{amsart}
\usepackage[T1]{fontenc}
\usepackage{amsmath,amssymb,mathtools,mathrsfs}
\usepackage{microtype}
\usepackage{booktabs,array,flafter,needspace}
\usepackage{tikz}
\usetikzlibrary{calc}
\usepackage[hidelinks]{hyperref}

\usepackage{enumitem}
\setlist[enumerate]{leftmargin=*,itemsep=3pt,topsep=5pt}
\numberwithin{equation}{section}
\newtheorem{theorem}{Theorem}[section]
\newtheorem{proposition}[theorem]{Proposition}
\newtheorem{lemma}[theorem]{Lemma}

\theoremstyle{definition}

\theoremstyle{remark}
\newtheorem{remark}[theorem]{Remark}
\DeclareMathOperator{\Isom}{Isom}
\DeclareMathOperator{\Comm}{Comm}

\newcommand{\R}{\mathbb R}
\newcommand{\C}{\mathbb C}
\newcommand{\Z}{\mathbb Z}
\newcommand{\HH}{\mathbb H}
\newcommand{\TT}{\mathbb T}
\newcommand{\cE}{\mathcal E}
\newcommand{\cD}{\mathcal D}

\newcommand{\cA}{\mathcal A}
\newcommand{\cS}{\mathcal S}
\newcommand{\cP}{\mathcal P}
\newcommand{\ip}[2]{\langle #1,#2\rangle}
\newcommand{\norm}[1]{\lVert #1\rVert}
\title[Hidden symmetries and echolocation]{Hidden symmetries and echolocation}

\author{Shilei Fan}
\address{School of Mathematics and Statistics, Central China Normal University, Wuhan 430079, PR China}
\email{slfan@mail.ccnu.edu.cn}
\author{Shi-Lei Kong}
\address{School of Mathematics, Sichuan University, Chengdu 610065, PR China}
\email{slkong@scu.edu.cn}

\author{Yakun Xi}
\address{School of Mathematical Sciences, Zhejiang University, Hangzhou 310027, PR China}
\email{yakunxi@zju.edu.cn}
\begin{document}
\begin{abstract}
We study the echolocation problem: whether local spectral data determine a point up to isometry. We characterize cospectrality by the existence of a spectral symmetry exchanging the two point evaluations and introduce a pointwise analogue of Sunada's method for constructing such symmetries. This yields cospectral, non-similar points on closed hyperbolic surfaces of every genus at least three and on round lens spaces. We also obtain Neumann examples on a scalene triangle and on a planar
slit domain, and
reinterpret a known Dirichlet example of Dixon, Miah, Safronov,
Snyder, Vakaryuk, and Wyman through the same framework.
\end{abstract}

\maketitle
\section{Introduction}
Let $(M,{\mathbf g})$ be a compact connected Riemannian manifold. If $M$ has
boundary, we impose either Dirichlet or Neumann boundary conditions on
the Laplacian. Write
\[
 -\Delta_{\mathbf g} e_j=\lambda_j^2e_j,
 \qquad 0\leq\lambda_1\leq\lambda_2\leq\cdots,
\]
where $\{e_j\}$ is a real orthonormal basis of $L^2(M)$. The pointwise
Weyl counting function at $x\in M$ is
\begin{equation}\label{eq:counting}
 N_x(\lambda)=\sum_{\lambda_j\leq\lambda}|e_j(x)|^2.
\end{equation}
We call two points $x,y$ \emph{cospectral} if $N_x=N_y$, and
\emph{similar} if an isometry of $(M,{\mathbf g})$ takes $x$ to $y$. These notions
are independent of the choice of orthonormal eigenbasis.

Similar points are necessarily cospectral. Wyman and Xi \cite{WX23}
asked whether the converse holds, phrasing this as an echolocation
problem: can one ``hear'' one's location from the local spectral data
$N_x$, at least up to isometry? For the Dirichlet problem, the points
under consideration are taken in the interior, since $N_x$ vanishes
identically at every boundary point. For the Neumann problem, boundary
points may also be considered, with evaluations understood by
continuous extension. A related probabilistic location-detection problem, using exit-time data for Brownian motion in place of local spectral data, was recently studied by Kratz and Langford \cite{KL26}.

Wyman and Xi \cite[Theorem 1.3]{WX23} proved that, on every closed
manifold of dimension at least two, a residual set of smooth metrics
has the stronger property that $N_x=N_y$ implies $x=y$. Further
examples and variants of the echolocation problem, including results
for boxes and mixed boundary conditions, were studied by
Wang \cite{Wang24}. F. Wang, Wyman, and Xi \cite{WWX25} subsequently
proved that a closed surface on which all points are cospectral must
be homogeneous. In a related inverse problem, X. Wang, Wyman, and
Xi \cite{WWX26} showed that the collection of pointwise Weyl counting
functions determines the Riemannian metric when the Laplace spectrum
is simple. The graph analogue of the echolocation problem and its
relation to cospectral vertices were studied in \cite{KWX25}.

The echolocation problem nevertheless fails for special geometries.
Dixon, Miah, Safronov, Snyder, Vakaryuk, and Wyman \cite{DMS26}
gave the first known connected counterexample: for the Dirichlet
Laplacian on a $30$--$60$--$90$ triangle, they exhibited two interior
cospectral points that are not related by an isometry.

Our first result identifies the symmetry that replaces isometry in the
echolocation problem. Although equality of the local spectral data does
not in general imply that two points are related by an isometry, we
show that it is equivalent to the existence of an orthogonal symmetry
of the Laplacian that exchanges their point evaluations. Thus
cospectrality always detects a symmetry, but one acting on functions
rather than necessarily on the underlying manifold. Let $\cE$ denote
the space of finite linear combinations of Laplace eigenfunctions.

\begin{theorem}[Spectral symmetry]\label{thm:spectral-symmetry}
For points $x,y\in M$ as above, the following are equivalent:
\begin{enumerate}
\item The points $x$ and $y$ are cospectral.
\item There exists an orthogonal involution $S$ on real $L^2(M)$
preserving every Laplace eigenspace such that
\begin{equation}\label{eq:intro-exchange}
 (Su)(x)=u(y),\qquad (Su)(y)=u(x)
 \qquad\text{for every }u\in\mathcal E.
\end{equation}
\end{enumerate}
\end{theorem}

We call an operator $S$ satisfying condition (2) of
Theorem \ref{thm:spectral-symmetry} a
\emph{spectral symmetry exchanging $x$ and $y$}.
Such an operator is self-adjoint, satisfies $S^2=I$,
and commutes with the Laplacian under the prescribed
boundary condition. Importantly, $S$ acts on the
function space $L^2(M)$ and need not be induced by
any map $M\to M$, let alone an isometry.
Thus Theorem \ref{thm:spectral-symmetry}
characterizes cospectrality in terms of a symmetry
of the function space rather than necessarily a
geometric symmetry of $M$.

The characterization is elementary. On each Laplace eigenspace,
cospectrality says that the two vectors representing evaluation at
$x$ and $y$ have the same norm. An orthogonal reflection exchanges
them when they differ; otherwise we take the identity. The direct
sum over all eigenspaces gives $S$.
Conversely, a simultaneous orthonormal eigenbasis for $S$ and the
Laplacian gives equality of the squared values at the two points.
The corresponding characterization for graph vertices appears in
Godsil and Smith \cite[Lemma 10.1]{GS17}. We give the proof in
Section \ref{sec:operators}.

The abstract existence of a spectral symmetry is therefore not the
geometric issue: the construction above already uses the local
spectral data. The problem is to produce such an operator directly
from the geometry, without first computing the eigenfunctions, while
ensuring that no isometry exchanges the marked points.

We develop a general mechanism for doing this. 
Geometric symmetries of a larger space, such as a
Riemannian cover or an unfolded domain, can be used to construct spectral symmetries on the original space, even when those geometric symmetries do not descend to isometries of the original space. Using this mechanism, we construct examples on closed manifolds. To our knowledge, these are the first counterexamples to the echolocation problem on closed manifolds. For every genus $g\geq3$, we obtain a closed hyperbolic surface containing $g-2$ pairs of cospectral, non-similar points, and we obtain a family of analogous examples on round lens spaces.
We also construct a pair of cospectral boundary points for the
Neumann Laplacian on a scalene triangle, a pair of cospectral interior
points for the Neumann Laplacian on a planar slit domain. We then give a
signed-translation proof of the Dirichlet example of \cite{DMS26}.
These apparently different counterexamples are all instances of the
same spectral-symmetry principle.

The closed-manifold constructions are based on the following finite
averaging criterion.

\begin{theorem}[Pointwise Sunada principle]\label{thm:orbit}
Let a finite group $G$ act isometrically on a closed connected
Riemannian manifold $X$, and let $H\leq G$ act freely. Suppose that one
$G$-orbit is the disjoint union of two distinct $H$-orbits:
\begin{equation}\label{eq:two-orbits}
Gz=Hz\sqcup Hw.
\end{equation}
Then the points $Hz,Hw$ of $M=X/H$ are cospectral. Averaging over $G$
induces an orthogonal projection $P$ on $L^2(M)$ such that
\begin{equation}\label{eq:reflection-average}
S=2P-I
\end{equation}
is a spectral symmetry exchanging the two point evaluations.
\end{theorem}

The averaging operator descends to $X/H$ even when the individual
elements of $G$ do not. Theorem \ref{thm:orbit} may thus be viewed
as a pointwise analogue of Sunada's method \cite{Sun85}, although
it is not an instance of Sunada's criterion. Both constructions
use finite-group averaging. In Sunada's method, two almost-conjugate
subgroups $H_1,H_2\leq G$ have averaging projections with equal
trace on every finite-dimensional $G$-representation, so their
invariant subspaces have equal dimension. Applied to Laplace
eigenspaces on a common cover, this gives equal eigenvalue
multiplicities on $X/H_1$ and $X/H_2$.

Theorem \ref{thm:orbit} instead fixes a single subgroup $H$ and a single
quotient $X/H$. In place of almost conjugacy, its additional input is
the pointwise orbit decomposition \eqref{eq:two-orbits}. The projection
onto the $G$-invariant subspace then identifies the two marked
evaluations, and reflection across its range exchanges them. Thus the
same finite-averaging mechanism that in Sunada's method compares global
spectral multiplicities on two quotients here compares local spectral
data at two points of one quotient.
Our first application is to closed hyperbolic surfaces.

\begin{theorem}[Hyperbolic surfaces]\label{thm:surface}
For every integer $g\geq3$, there is a closed connected orientable
hyperbolic surface $M$ of genus $g$ containing $g-2$ pairs of distinct
points
\[
 x_1,y_1,\ldots,x_{g-2},y_{g-2}
\]
such that
\[
 N_{x_i}=N_{y_i}
 \qquad (1\leq i\leq g-2),
\]
but no isometry of $M$ carries $x_i$ to $y_i$.
\end{theorem}

The construction realizes the pointwise Sunada principle using the
dihedral group $D_4$ of order eight. We begin with a hyperbolic surface
$X$ of genus $2g-1$ carrying an isometric $D_4$-action, and quotient by
a freely acting subgroup $H\leq D_4$ of order two to obtain a surface
$M=X/H$ of genus $g$. Certain $D_4$-orbits on $X$ split into two
$H$-orbits. They descend to the pairs $x_i,y_i$ on $M$, and the
orbit-reflection principle gives their cospectrality. This orbit
structure persists throughout the family of hyperbolic realizations
carrying the prescribed $D_4$-action.
We then choose a generic hyperbolic realization
so that no isometry of $M$ exchanges any of
the marked pairs. 
There is also an explicit family in dimension three. We work on the
round $3$-sphere $S^3$, viewed in complex coordinates $(z_1,z_2)\in\mathbb C^2$. Let
$m\geq3$, let $q$ be relatively prime to $m$, and set
$ \zeta=e^{2\pi i/m}.$

Consider the isometry
\[
 \gamma(z_1,z_2)=(\zeta z_1,\zeta^q z_2).
\]
The cyclic group $\langle\gamma\rangle$ has order $m$ and acts freely
on $S^3$. We write
\[
 L(m;1,q)=S^3/\langle\gamma\rangle
\]
for the resulting lens space, equipped with the quotient of the round
metric.

\begin{theorem}[Lens spaces]\label{thm:lens}
Let $x,y\in L(m;1,q)$ be the images of
$(1,0)$ and $(0,1)$ 
respectively. Then $x$ and $y$ are cospectral. They are similar if and
only if
\begin{equation}\label{eq:lens-congruence}
 q^2\equiv1\quad\text{or}\quad q^2\equiv-1\pmod m.
\end{equation}
In particular, $L(7;1,2)$ contains a pair of cospectral, non-similar
points.
\end{theorem}
Here the orbit-reflection principle has a simple geometric interpretation.
On $S^3$, a larger group of symmetries moves the points lying over $x$
and $y$ through a single orbit. Averaging over these symmetries therefore
treats the two points in the same way, and the orbit-reflection principle
gives their cospectrality. These symmetries need not descend to actual isometries of the lens space.

We next turn to planar examples. Let $T$ be the open $30$--$60$--$90$
triangle with vertices
\[
 A=(0,0),\qquad
 B=\left(\frac{\sqrt3}{2},0\right),\qquad
 C=\left(\frac{\sqrt3}{2},\frac12\right).
\]
Since $T$ is scalene, its isometry group is trivial.

\begin{theorem}[Neumann case]\label{thm:neumann}
For the Neumann Laplacian on $T$, the boundary points
\[
 P_N=\left(\frac{\sqrt3}{4},0\right),\qquad
 Q_N=\left(\frac{3\sqrt3}{8},\frac38\right)
\]
are cospectral and non-similar.
\end{theorem}
\begin{figure}[htbp]
\centering
\begin{tikzpicture}[x=8cm,y=8cm,font=\small]
\coordinate (A) at (0,0);
\coordinate (B) at ({sqrt(3)/2},0);
\coordinate (C) at ({sqrt(3)/2},1/2);

\coordinate (pN) at ({sqrt(3)/4},0);
\coordinate (qN) at ({3*sqrt(3)/8},3/8);
\coordinate (pD) at ({7*sqrt(3)/16},3/16);
\coordinate (qD) at ({5*sqrt(3)/16},3/16);

\draw[thick] (A)--(B)--(C)--cycle;

\fill (pN) circle (1.5pt);
\fill (qN) circle (1.5pt);
\fill (pD) circle (1.5pt);
\fill (qD) circle (1.5pt);

\node[below left] at (A) {$A$};
\node[below right] at (B) {$B$};
\node[above right] at (C) {$C$};

\node[below] at (pN) {$P_N$};
\node[above left] at (qN) {$Q_N$};
\node[below right] at (pD) {$P_D$};
\node[below left] at (qD) {$Q_D$};
\end{tikzpicture}
\caption{The triangle $T$ with the Neumann pair $P_N,Q_N$ on the
boundary and the Dirichlet pair $P_D,Q_D$ in the interior.}
\label{fig:triangle-points}
\end{figure}
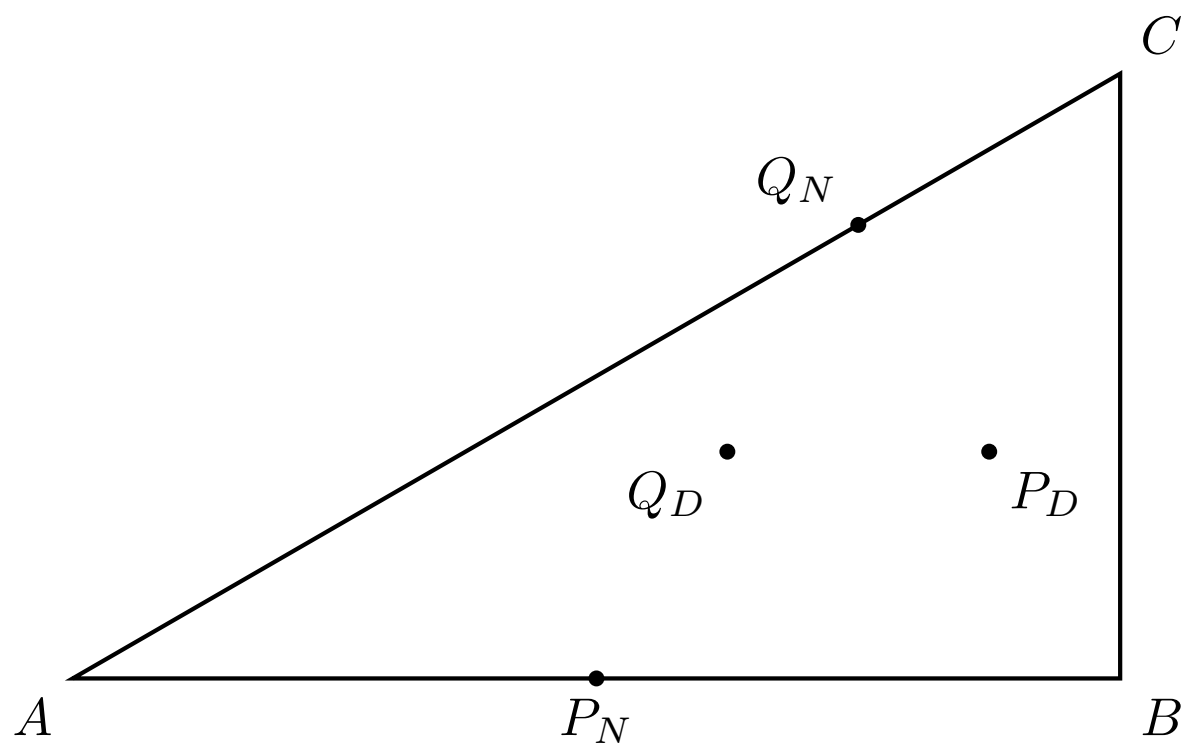
The Neumann example is obtained by unfolding the triangle to a flat
torus, where averaging over finitely many symmetries produces a spectral
symmetry exchanging the two marked points. Although
Theorem \ref{thm:orbit} does not apply directly, the construction is
driven by the same averaging-reflection mechanism.

For the Dirichlet problem on $T$, Dixon, Miah, Safronov, Snyder,
Vakaryuk, and Wyman \cite{DMS26} found the cospectral pair
\begin{equation}\label{eq:dirichlet-points}
 P_D=\left(\frac{7\sqrt3}{16},\frac3{16}\right),\qquad
 Q_D=\left(\frac{5\sqrt3}{16},\frac3{16}\right).
\end{equation}
We reinterpret their example using odd reflection to the same flat
torus. A signed combination of translations induces a self-adjoint
operator that preserves the Dirichlet eigenspaces and exchanges
the two marked evaluations. A self-adjoint exchange criterion
then gives cospectrality.

We further show that both triangle pairs are strongly cospectral
and exhibit perfect Schr\"odinger state transfer at explicit times
(Proposition \ref{prop:triangle-strong-cospectrality}). These notions naturally extend their graph-theoretic
counterparts \cite{GS17} to the continuous setting.

There is also a simple planar Neumann example with both marked points in the
interior. The price is that the domain has a slit.

Let
\[
 \Omega=
 \bigl((-1,1)\times(0,3)\bigr)
 \setminus
 \left(\{0\}\times\left[0,\frac12\right]\right),
\]
with Neumann boundary conditions on the outer rectangle and on both
faces of the slit.

\begin{theorem}[Slit domain]\label{thm:slit}
For the Neumann Laplacian on $\Omega$, the interior points
\[
 p=(0,1),\qquad q=(0,2)
\]
are cospectral and non-similar.
\end{theorem}

\begin{remark}[Relation to transplantation]
The closed examples apply Theorem \ref{thm:orbit} directly.
The triangle examples use related averaging and signed-translation
operators, while the slit example uses an even--odd decomposition.

A related precedent is the homophonic transplantation of
Buser, Conway, Doyle, and Semmler \cite{BCDS94}, which matches
point evaluations at distinguished points in two \emph{different} drums.
Our constructions produce analogous matching evaluations within
a single spectral problem.
\end{remark}

\subsection*{Structure of the paper}

Section \ref{sec:operators} proves
Theorem \ref{thm:spectral-symmetry} and records the self-adjoint
exchange criterion. Section \ref{sec:averaging} proves the pointwise
Sunada principle. We treat lens spaces and hyperbolic surfaces in
Sections \ref{sec:lens} and \ref{sec:surface}, and the Neumann and
Dirichlet triangle constructions in Section \ref{sec:triangles}.
Finally, Section \ref{sec:slit} gives the Neumann slit-domain example. 

\subsection*{Statement on AI use}
This project began with our use of ChatGPT 5.6 sol to search for a negative example to the echolocation problem. ChatGPT 5.6 sol suggested the lens-space example. We then recognized that a Sunada-type principle was at work and used ChatGPT 5.6 sol to explore analogous examples on hyperbolic surfaces and  the Neumann problem on a triangle. We subsequently formulated and summarized the ideas that became Theorems \ref{thm:spectral-symmetry} and \ref{thm:orbit} and developed them into the present paper. ChatGPT was also used to polish and proofread the manuscript. 
The authors view this work as an example of LLM-assisted
``experimental mathematics'', in which an LLM suggests examples, while the authors identify the underlying structures and formulate general principles.

\subsection*{Acknowledgment}  This project is supported by the National Key Research and Development Program of China No. 2022YFA1007200. Y. X. is supported by Zhejiang Provincial Natural Science Foundation of China under
Grant No. LR25A010001, and NSFC Grant No.
12571107.

\section{Operators exchanging point evaluations}\label{sec:operators}

We first prove Theorem \ref{thm:spectral-symmetry}. We then record a
weaker exchange criterion, which will be useful later.

For each distinct frequency $\mu$, let
\[
 E_\mu=\{u:-\Delta_{\mathbf g} u=\mu^2u\},
\]
with the prescribed boundary condition, and let $\Pi_\mu$ denote the
orthogonal projection onto $E_\mu$, as well as its Schwartz kernel.
Then
\[
 \Pi_\mu(x,x)=\sum_{\lambda_j=\mu}|e_j(x)|^2,
\]
so $x$ and $y$ are cospectral if and only if
\[
 \Pi_\mu(x,x)=\Pi_\mu(y,y)
 \qquad\text{for every }\mu.
\]
Although point evaluation is not in general bounded on $L^2(M)$, it is
well defined on each finite-dimensional eigenspace. All evaluation identities below are understood on $\cE$, the space of
finite linear combinations of Laplace eigenfunctions. All inner products and norms below are those of real $L^2(M)$.

\begin{proof}[Proof of Theorem \ref{thm:spectral-symmetry}]
Suppose first that $x$ and $y$ are cospectral. For each $\mu$, let
$a_\mu,b_\mu\in E_\mu$ be the vectors representing evaluation at $x$
and $y$, that is
\[
 u(x)=\ip{u}{a_\mu},\qquad
 u(y)=\ip{u}{b_\mu},
 \qquad u\in E_\mu.
\]
Explicitly,
\[
 a_\mu=\sum_{\lambda_j=\mu}e_j(x)e_j,\qquad
 b_\mu=\sum_{\lambda_j=\mu}e_j(y)e_j.
\]
Cospectrality is exactly the statement that
\begin{equation}\label{eq:evaluation-norms}
 \norm{a_\mu}^2=\Pi_\mu(x,x)
 =\Pi_\mu(y,y)=\norm{b_\mu}^2.
\end{equation}

If $a_\mu=b_\mu$, set $S_\mu=I$ on $E_\mu$. Otherwise, put
$d_\mu=a_\mu-b_\mu$ and define
\begin{equation}\label{eq:eigenspace-reflection}
 S_\mu u
 =u-2\frac{\ip{u}{d_\mu}}{\norm{d_\mu}^2}d_\mu.
\end{equation}
This is the orthogonal reflection across $d_\mu^\perp$. Since
$\norm{a_\mu}=\norm{b_\mu}$,
\[
 \ip{a_\mu}{d_\mu}
 =\frac12\norm{d_\mu}^2,
\]
and hence
\[
 S_\mu a_\mu=b_\mu,\qquad
 S_\mu b_\mu=a_\mu.
\]

Now set
\[
 S=\bigoplus_\mu S_\mu.
\]
Then $S$ is an orthogonal involution on real $L^2(M)$ preserving every
Laplace eigenspace. For $u=\sum_\mu u_\mu\in\cE$, self-adjointness of
each $S_\mu$ gives
\[
 (Su)(x)
 =\sum_\mu\ip{S_\mu u_\mu}{a_\mu}
 =\sum_\mu\ip{u_\mu}{S_\mu a_\mu}
 =\sum_\mu\ip{u_\mu}{b_\mu}
 =u(y).
\]
Similarly,
\[
 (Su)(y)=u(x).
\]

Conversely, suppose that $S$ has the stated properties. Since $S$ is an
orthogonal involution, it is self-adjoint. On each $E_\mu$, choose an
orthonormal basis of eigenvectors of $S$. Thus every vector $u$ in the
resulting Laplace eigenbasis satisfies
\[
 Su=\varepsilon u,
 \qquad \varepsilon\in\{1,-1\}.
\]
The exchange identity gives
\[
 u(y)=(Su)(x)=\varepsilon u(x),
\]
and hence $|u(x)|^2=|u(y)|^2$. Summing over an orthonormal basis of
each $E_\mu$ yields
\[
 \Pi_\mu(x,x)=\Pi_\mu(y,y)
\]
for every $\mu$. Thus $x$ and $y$ are cospectral.
\end{proof}

For some constructions it is easier to produce a self-adjoint exchange
operator before normalizing it to an involution. We therefore record
the following weaker criterion.

\begin{lemma}[Self-adjoint exchange criterion]\label{lem:exchange}
Suppose that a bounded self-adjoint operator $A$ on real $L^2(M)$
preserves every $E_\mu$ and satisfies
\[
 (Au)(x)=u(y),\qquad
 (Au)(y)=u(x),
 \qquad u\in\cE.
\]
Then $x$ and $y$ are cospectral.
\end{lemma}

\begin{proof}
On $E_\mu$, the exchange identities and self-adjointness give
\[
 Aa_\mu=b_\mu,\qquad
 Ab_\mu=a_\mu.
\]
Therefore
\[
 \norm{a_\mu}^2
 =\ip{Ab_\mu}{a_\mu}
 =\ip{b_\mu}{Aa_\mu}
 =\norm{b_\mu}^2.
\]
Thus
\[
 \Pi_\mu(x,x)=\Pi_\mu(y,y)
\]
for every $\mu$, and hence $x$ and $y$ are cospectral.
\end{proof}

\section{Averaging on a finite quotient}\label{sec:averaging}

We prove Theorem \ref{thm:orbit}. For $F\in L^2(X)$, define the
$G$-average
\begin{equation}\label{eq:group-average}
 P_GF=\frac1{|G|}\sum_{g\in G}F\circ g.
\end{equation}
Since each pullback $F\mapsto F\circ g$ is unitary and commutes with
the Laplacian, the group average satisfies
\[
 P_G^*=P_G,\qquad P_G^2=P_G.
\]
Hence $P_G$ is the orthogonal projection onto the $G$-invariant
functions,
\[
 L^2(X)^G
 :=
 \{F\in L^2(X):F\circ g=F\text{ for every }g\in G\}.
\]

Now let $\pi:X\to M=X/H$ be the quotient map. Functions on $M$ lift
precisely to $H$-invariant functions on $X$. More precisely,
\[
 Uu=|H|^{-1/2}u\circ\pi
\]
defines a unitary map
\[
 U:L^2(M)\longrightarrow L^2(X)^H.
\]

Since $H\leq G$, every $G$-invariant function is $H$-invariant, so
\[
 L^2(X)^G\subset L^2(X)^H.
\]
Thus $P_G$ restricts to an orthogonal projection on the $H$-invariant
subspace and, through the identification $U$, induces an orthogonal
projection on $L^2(M)$:
\[
 P=U^{-1}P_GU.
\]
Finally, set
\[
 S=2P-I.
\]
Since $P$ is an orthogonal projection, $S$ is the reflection across
the range of $P$. In particular, $S$ is an orthogonal involution.
Moreover, because $P_G$ commutes with the Laplacian, $S$ preserves
every Laplace eigenspace on $M$.

It remains to check the two marked evaluations. Let
\[
 x=Hz,\qquad y=Hw,
\]
and write
\[
 a=u(x),\qquad b=u(y).
\]
The lift $F=u\circ\pi$ is constant with value $a$ on $Hz$ and with
value $b$ on $Hw$. By hypothesis,
\[
 Gz=Hz\sqcup Hw,
\]
and, since the $H$-action is free, the two $H$-orbits have the same
cardinality. Averaging over $G$ therefore gives
\[
 (P_GF)(z)=(P_GF)(w)=\frac{a+b}{2}.
\]
Equivalently, downstairs,
\[
 Pu(x)=Pu(y)=\frac{a+b}{2}.
\]
Hence
\[
 (Su)(x)=2Pu(x)-u(x)=b=u(y),
\]
and similarly $(Su)(y)=u(x)$. Thus $S$ exchanges the two point
evaluations, and Theorem \ref{thm:spectral-symmetry} gives the
cospectrality of $x$ and $y$.

\begin{remark}
An element $g\in G$ need not induce a well-defined map on the quotient
$X/H$. Indeed, a point of $X/H$ is an $H$-orbit, so for $g$ to descend
it must send $H$-orbits to $H$-orbits. Equivalently, one needs
\[
 gHg^{-1}=H.
\]
Thus an element of $G$ may carry $z$ to $w$ upstairs without giving an
isometry of $X/H$ carrying $Hz$ to $Hw$.

The full average behaves differently. If $F$ is $H$-invariant, then
$P_GF$ is $G$-invariant, and hence in particular $H$-invariant.
Therefore the averaging operator descends to $X/H$ even when the
individual symmetries being averaged do not.

This distinction is important in our constructions. The spectral
symmetry can descend from the averaging construction without coming
from an actual isometry of the quotient. Proving that the marked points
are non-similar then requires ruling out any other accidental isometries.
\end{remark}

\section{Lens spaces}\label{sec:lens}

We first use Theorem \ref{thm:orbit} to prove Theorem \ref{thm:lens}. Let
\[
 H=\langle\gamma\rangle,
 \qquad
 \gamma(z_1,z_2)=(\zeta z_1,\zeta^qz_2).
\]
The action of $H$ on $S^3$ is free. Indeed, if $\gamma^k$ fixes a
point, then at least one coordinate is nonzero. If the first coordinate
is nonzero, then $k\equiv0\pmod m$; if the second is nonzero, then
$qk\equiv0\pmod m$, and again $k\equiv0\pmod m$ since
$\gcd(m,q)=1$.

Consider the finite group $G$ generated by the independent rotations
\[
 (z_1,z_2)\longmapsto(\zeta z_1,z_2),
 \qquad
 (z_1,z_2)\longmapsto(z_1,\zeta z_2),
\]
together with the coordinate interchange
\[
 (z_1,z_2)\longmapsto(z_2,z_1).
\]
Let $v=(1,0)$ and $w=(0,1)$. 
The $G$-orbit of $v$ consists of the two coordinate circles sampled at
the $m$th roots of unity:
\[
 Gv
 =
 \{(\zeta^j,0):0\leq j<m\}
 \sqcup
 \{(0,\zeta^j):0\leq j<m\}.
\]
The first set is exactly the $H$-orbit of $v$. Similarly
\[
 Hw=\{(0,\zeta^j):0\leq j<m\}.
\]
Therefore
\[
 Gv=Hv\sqcup Hw.
\]
Theorem \ref{thm:orbit} now shows that the images $x$ and $y$ of $v$
and $w$ in $L(m;1,q)=S^3/H$ are cospectral.

It remains to determine when $x$ and $y$ are related by an isometry.
Suppose that an isometry of the lens space carries $x$ to $y$. It
lifts to an isometry $A\in O(4)$ of $S^3$ that normalizes the deck group $H$.
After composing with a deck transformation, we may assume that
$Av=w$. Since $A$ normalizes the cyclic group $H=\langle\gamma\rangle$,
there is an integer $r$, relatively prime to $m$, such that
\[
 A\gamma A^{-1}=\gamma^r.
\]

The $H$-orbits of $v$ and $w$ lie in the two coordinate planes
\[
 V_1=\C\times\{0\},
 \qquad
 V_2=\{0\}\times\C,
\]
and each orbit spans its coordinate plane over $\R$. Since $A$ carries
the $H$-orbit of $v$ to the $H$-orbit of $w$, it follows that
\[
 A(V_1)=V_2,\qquad A(V_2)=V_1.
\]

On $V_1$, the generator $\gamma$ is rotation through angle $2\pi/m$,
while on $V_2$ it is rotation through angle $2\pi q/m$. Two rotations
of a real plane are orthogonally conjugate exactly when their angles
agree up to sign. Thus
\begin{equation}\label{eq:rotation-congruences}
 1\equiv\varepsilon qr\pmod m,
 \qquad
 q\equiv\eta r\pmod m,
 \qquad
 \varepsilon,\eta\in\{1,-1\}.
\end{equation}
Therefore
\[
 q^2\equiv\pm1\pmod m.
\]

Conversely, suppose first that $q^2\equiv1\pmod m$. Then it is routine to check that the coordinate
interchange
$(z_1,z_2)\longmapsto(z_2,z_1)
$
normalizes $H$ and exchanges $v$ and $w$. If instead
$q^2\equiv-1\pmod m$, the map
$(z_1,z_2)\longmapsto(\overline{z_2},z_1)
$
normalizes $H$ and also exchanges $v$ and $w$. In either case the map
descends to an isometry of the lens space carrying $x$ to $y$.

Thus $x$ and $y$ are similar exactly when
\[
 q^2\equiv\pm1\pmod m,
\]
which proves Theorem \ref{thm:lens}.

\section{Hyperbolic surfaces}\label{sec:surface}

We now turn to the hyperbolic surface construction, where the
orbit-reflection principle requires a more elaborate covering space.
Throughout this section,
\[
 D_4=\langle r,s:r^4=s^2=1,\ srs=r^{-1}\rangle
\]
is the dihedral group of order eight. The generators are abstract:
the action constructed below is orientation preserving, so $s$ should
not be interpreted as a reflection of the hyperbolic plane.

\subsection{The covering surface and the marked points}

Fix $g\geq3$. Choose for now a cocompact orientation-preserving
Fuchsian group $\Gamma$ of signature
\[
 (0;4,4,\underbrace{2,\ldots,2}_{g\text{ times}}).
\]
Thus $\HH^2/\Gamma$ is a sphere with two cone points of order four and
$g$ cone points of order two. We may choose elliptic generators
$a_1,\ldots,a_{g+2}$ with presentation
\begin{equation}\label{eq:orbifold-presentation}
 \Gamma=
 \langle a_1,\ldots,a_{g+2}:
 a_1^4=a_2^4=a_3^2=\cdots=a_{g+2}^2=1,\quad
 a_1\cdots a_{g+2}=1\rangle.
\end{equation}

Define a homomorphism
\[
 \phi:\Gamma\longrightarrow D_4
\]
by
\begin{equation}\label{eq:monodromy}
\begin{aligned}
 \phi(a_1)&=r,&
 \phi(a_2)&=r^{(-1)^{g+1}},\\
 \phi(a_3)&=\phi(a_4)=s,&
 \phi(a_j)&=r^2\qquad(5\leq j\leq g+2).
\end{aligned}
\end{equation}
Indeed, the product of the images is
\[
 r^{1+(-1)^{g+1}+2(g-2)}=1,
\]
since the exponent is divisible by four. Thus the defining relation is
preserved. Moreover, $r$ and $s$ occur among the images, so $\phi$ is
surjective. For $g=3$, the assignment is simply
\[
 (a_1,a_2,a_3,a_4,a_5)\longmapsto(r,r,s,s,r^2).
\]

Let
$
 X=\HH^2/\ker\phi.
$
Each elliptic generator of $\Gamma$ has the same order as its image
under $\phi$. Since every nontrivial finite-order element of $\Gamma$
is conjugate to a nontrivial power of an elliptic generator, no such
element lies in $\ker\phi$. Hence $\ker\phi$ is torsion free, so $X$ is a smooth
closed hyperbolic surface. Since
\[
 \Gamma/\ker\phi\cong D_4,
\]
the group $D_4$ acts on $X$ by orientation-preserving isometries.

Now let
\[
 h=rs,\qquad H=\langle h\rangle.
\]
The only nontrivial stabilizers in the $D_4$-action on $X$ are
conjugate to the cyclic groups generated by $r$, $s$, or $r^2$.
The involutions conjugate to $s$ are $s$ and $r^2s$, whereas
$h=rs$ belongs to the other conjugacy class
 $\{rs,r^3s\}.$
Thus no nontrivial element of $H$ fixes a point of $X$. In other words,
$H$ acts freely, and
 $M=X/H$
is a smooth closed orientable hyperbolic surface.

We next compute its genus. The orbifold Euler characteristic of
$\HH^2/\Gamma$ is
\[
 \chi_{\mathrm{orb}}(\HH^2/\Gamma)
 =
 2-2\left(1-\frac14\right)
 -g\left(1-\frac12\right)
 =
 \frac{1-g}{2}.
\]
Since $[\Gamma:\ker\phi]=8$, we have
\[
 \chi(X)=8\,\chi_{\mathrm{orb}}(\HH^2/\Gamma)=4-4g.
\]
The quotient $X\to M=X/H$ has degree two, so
\[
 \chi(M)=\frac{\chi(X)}2=2-2g.
\]
Thus $X$ has genus $2g-1$ and $M$ has genus $g$.

For each $i=1,\ldots,g-2$, choose a point $p_i\in X$ lying over the
cone point corresponding to $a_{i+4}$. Its stabilizer in $D_4$ is
exactly $\langle r^2\rangle$. Define
\[
 x_i=Hp_i,\qquad y_i=Hrp_i.
\]
Since $h=rs$, the four-point $D_4$-orbit of $p_i$ splits as
\begin{equation}\label{eq:surface-orbit}
 D_4p_i
 =
 \{p_i,hp_i\}
 \sqcup
 \{rp_i,hrp_i\}.
\end{equation}
These are two distinct $H$-orbits. Hence
Theorem \ref{thm:orbit} gives
\[
 N_{x_i}=N_{y_i}
 \qquad(1\leq i\leq g-2).
\]
Points arising from different values of $i$ lie over different cone
points of $\HH^2/\Gamma$, so all $2g-4$ marked points are distinct.

We have therefore constructed the required cospectral pairs for every
hyperbolic realization of the prescribed orbifold. It remains to
choose the realization so that no isometry of $M$ exchanges any
$x_i$ with $y_i$.

\subsection{Choosing the hyperbolic structure}

It remains to choose the hyperbolic realization so that the marked
points are not exchanged by accidental isometries. For this purpose we
choose $\Gamma$ with no hidden symmetries beyond those already contained
in $\Gamma$. Recall that the commensurator of $\Gamma$ in the full
isometry group of $\HH^2$ is
\[
 \Comm(\Gamma)
 =
 \{c\in\Isom(\HH^2):
 \Gamma\cap c\Gamma c^{-1}
 \text{ has finite index in both groups}\}.
\]

\begin{lemma}\label{lem:generic}
For every $g\geq3$, there is a cocompact Fuchsian group $\Gamma$ of
signature
\[
 (0;4,4,\underbrace{2,\ldots,2}_{g\text{ times}})
\]
such that
\begin{equation}\label{eq:commensurator}
 \Comm(\Gamma)=\Gamma.
\end{equation}
In fact, such realizations form a full-measure subset of the
$(g-1)$-dimensional complex deformation space of labeled orbifolds of
this signature.
\end{lemma}

\begin{proof}
A labeled orbifold of this signature is topologically a sphere with $g+2$ marked cone 
points, two of order four and $g$ of order two. After fixing three
marked points at $0,1,\infty$, the remaining configurations form an
open subset of $\C^{g-1}$. We show that the realizations for which
$\Comm(\Gamma)$ is larger than $\Gamma$ lie in a measure-zero subset.

We first exclude realizations for which $\Gamma$ is contained in a
proper discrete orientation-preserving overgroup $\Gamma'$. Suppose
$[\Gamma':\Gamma]=d\geq2$ and that $\Gamma'$ has signature
$(b;m_1,\ldots,m_r)$. Multiplicativity of orbifold area gives
\begin{equation}\label{eq:supergroup-area}
 2b-2+\sum_{j=1}^r\left(1-\frac1{m_j}\right)
 =
 \frac{g-1}{2d}
 \leq
 \frac{g-1}{4}.
\end{equation}
Because $\Gamma'$ contains an element of order four, at least one
$m_j$ is divisible by four. Hence the sum on the left is at least
$r/2+1/4$, and therefore $8b+2r\leq g+6$. It follows that the
deformation space of the quotient orbifold has complex dimension
\[
 3b-3+r\leq \frac g2-b<g-1.
\]
Thus every possible proper overgroup occurs in a family of strictly
smaller dimension. For fixed degree, quotient signature, and
topological covering type, the structure upstairs is obtained by
pulling back the structure downstairs, and there are only countably
many such covering types. Consequently, the realizations admitting a
proper discrete overgroup form a measure-zero subset.

We next exclude arithmetic realizations. By Borel's finiteness
theorem \cite[Theorem 8.2]{Bor81}, there are only finitely many
conjugacy classes of arithmetic Fuchsian groups with bounded covolume.
Since the signature is fixed, so is the hyperbolic area, and therefore
only finitely many arithmetic conjugacy classes occur in our deformation
space. In particular, they form a measure-zero subset.

Now choose a non-arithmetic realization outside the first exceptional
set. By Margulis's commensurator theorem \cite[Theorem 1]{Mar91},
the orientation-preserving commensurator
$\Comm_{\Isom^+(\HH^2)}(\Gamma)$ is discrete. Since it contains
$\Gamma$, it is a discrete orientation-preserving overgroup of
$\Gamma$. The first step rules out any proper such overgroup, and hence
\[
 \Comm_{\Isom^+(\HH^2)}(\Gamma)=\Gamma.
\]

Finally, we exclude realizations admitting an orientation-reversing
isometry. Such an isometry must permute the cone points while preserving
their orders. For each of the finitely many possible permutations,
requiring such a symmetry imposes a nontrivial condition on their
positions, so the corresponding configurations form a measure-zero
subset.

Choose a realization outside all three exceptional sets. We already
know that the orientation-preserving commensurator is $\Gamma$. Suppose
that the full commensurator contained an orientation-reversing element
$c$. Since the commensurator is a group, we have
\[
 c\Gamma c^{-1}\subset
 \Comm_{\Isom^+(\HH^2)}(\Gamma)=\Gamma.
\]
Thus $c$ normalizes $\Gamma$ and descends to an
orientation-reversing isometry of $\HH^2/\Gamma$, contrary to our
choice. Hence $\Comm(\Gamma)=\Gamma$.

\end{proof}

\subsection{The isometry group of $M$}

Choose $\Gamma$ as in Lemma \ref{lem:generic}, and let
\[
 \Lambda=\phi^{-1}(H).
\]
Since $M=X/H$ and $X=\HH^2/\ker\phi$, we may identify
\[
 M=\HH^2/\Lambda.
\]
Every isometry of $M$ lifts to an isometry of $\HH^2$ normalizing
$\Lambda$, so
\[
 \Isom(M)=N_{\Isom(\HH^2)}(\Lambda)/\Lambda.
\]
Here $N_{\Isom(\HH^2)}(\Lambda)$ is the normalizer of $\Lambda$ in
$\Isom(\HH^2)$.

Since $\Lambda$ has finite index in $\Gamma$, any isometry normalizing
$\Lambda$ also commensurates $\Gamma$. By our choice of $\Gamma$,
$\Comm(\Gamma)=\Gamma$, and therefore
\[
 \Isom(M)=N_\Gamma(\Lambda)/\Lambda.
\]
Now $\Lambda=\phi^{-1}(H)$ and $\ker\phi\subset\Lambda$, so the quotient
map $\phi:\Gamma\to D_4$ gives
\begin{equation}\label{eq:full-isometry}
 \Isom(M)\cong N_{D_4}(H)/H.
\end{equation}

For $H=\langle rs\rangle$, a direct calculation gives
\[
 N_{D_4}(H)=\{1,rs,r^2,r^3s\}.
\]
Thus $\Isom(M)$ has order two, with its nontrivial element represented
by $r^2$.

For each $i$, the stabilizer of $p_i$ in $D_4$ is
$\langle r^2\rangle$, so $r^2p_i=p_i$. Since $r^2$ commutes with $r$,
it also fixes $rp_i$. The induced isometry of $M$ therefore fixes
$x_i=Hp_i$ and $y_i=Hrp_i$ separately. In particular, no isometry of
$M$ carries $x_i$ to $y_i$.
Together with $N_{x_i}=N_{y_i}$ for $1\leq i\leq g-2$, this gives
$g-2$ pairs of cospectral, non-similar points and completes the proof
of Theorem \ref{thm:surface}.

\section{A triangle}\label{sec:triangles}

Both triangle examples can be described on the same unfolded flat
torus. Recall that $T$ is the open $30$--$60$--$90$ triangle with
vertices
\[
 A=(0,0),\qquad
 B=\left(\frac{\sqrt3}{2},0\right),\qquad
 C=\left(\frac{\sqrt3}{2},\frac12\right).
\]
Reflections in its sides tile the plane. The translation subgroup
of the resulting reflection group is the lattice
\[
 \Lambda=\Z t_1+\Z t_2,\qquad
 t_1=(\sqrt3,0),\qquad
 t_2=\left(\frac{\sqrt3}{2},\frac32\right).
\]

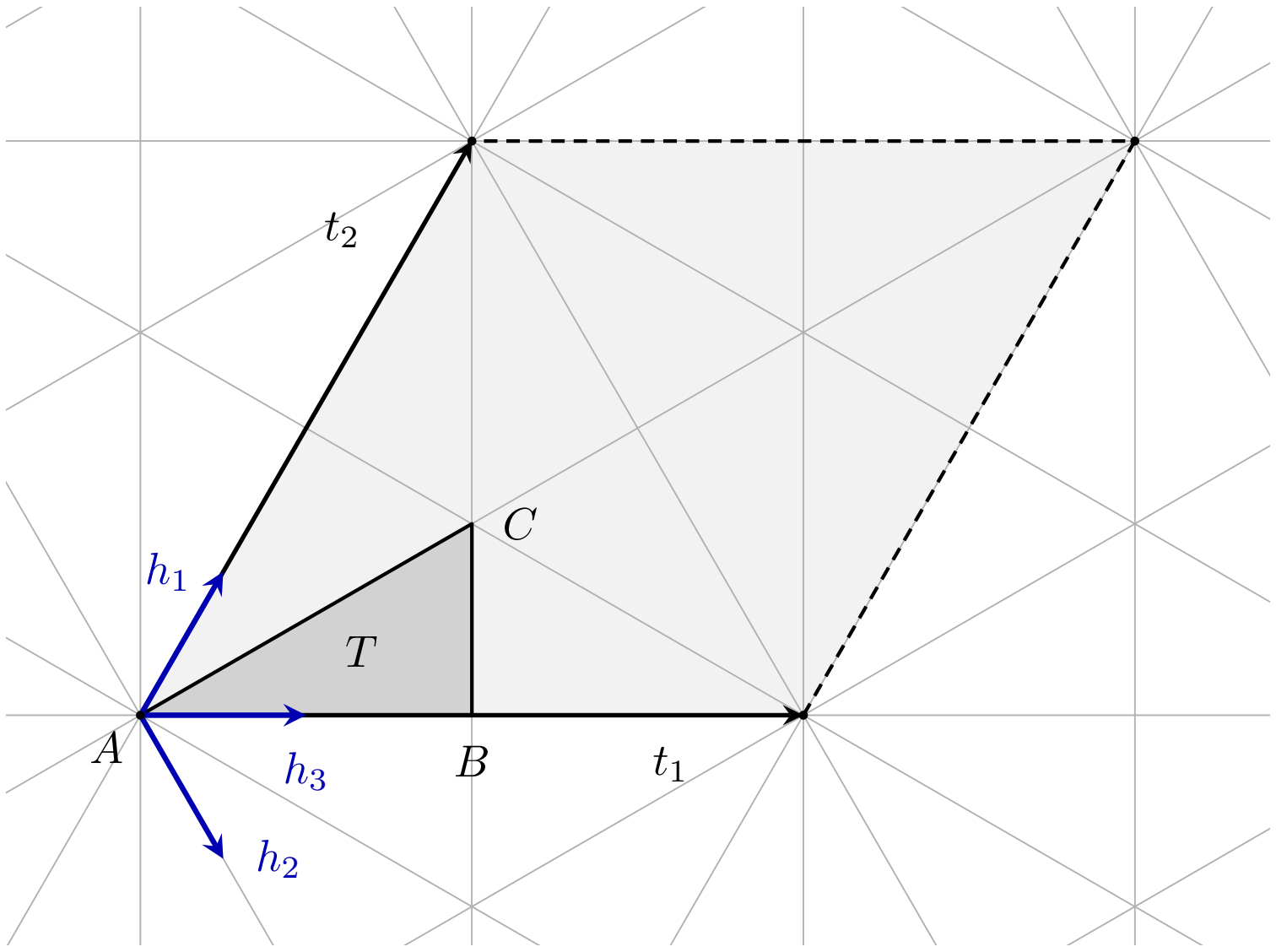
\begin{figure}[htbp]
\centering
\begin{tikzpicture}[x=3cm,y=3cm,font=\small,line join=round]
\coordinate (A) at (0,0);
\coordinate (B) at ({sqrt(3)/2},0);
\coordinate (C) at ({sqrt(3)/2},.5);
\coordinate (L1) at ({sqrt(3)},0);
\coordinate (L2) at ({sqrt(3)/2},1.5);
\coordinate (L12) at ({3*sqrt(3)/2},1.5);
\coordinate (H1) at ({sqrt(3)/8},{3/8});
\coordinate (H2) at ({sqrt(3)/8},{-3/8});
\coordinate (H3) at ({sqrt(3)/4},0);

\begin{scope}
  \clip (-.35,-.60) rectangle (2.95,1.85);
  \fill[gray!10] (A)--(L1)--(L12)--(L2)--cycle;
  \foreach \tilei in {-1,0,1,2} {
    \foreach \tilej in {-1,0,1,2} {
      \begin{scope}[
        shift={({sqrt(3)*(\tilei+\tilej/2)},{1.5*\tilej})}]
        \foreach \tileangle in {0,60,120,180,240,300} {
          \begin{scope}[rotate=\tileangle]
            \draw[gray!60,line width=.3pt]
              (0,0)--({sqrt(3)/2},0)--
              ({sqrt(3)/2},.5)--cycle;
            \draw[gray!60,line width=.3pt]
              (0,0)--({sqrt(3)/2},0)--
              ({sqrt(3)/2},-.5)--cycle;
          \end{scope}
        }
      \end{scope}
    }
  }
\end{scope}

\filldraw[fill=gray!35,draw=black,line width=.8pt]
  (A)--(B)--(C)--cycle;
\node at ({sqrt(3)/3},{1/6}) {$T$};
\draw[densely dashed,line width=.8pt] (L1)--(L12)--(L2);
\draw[->,>=stealth,line width=1pt]
  (A)--(L1) node[pos=.80,below=3pt] {$t_1$};
\draw[->,>=stealth,line width=1pt]
  (A)--(L2) node[pos=.76,above left=3pt] {$t_2$};

\begin{scope}[draw=blue!70!black,text=blue!70!black]
  \draw[->,>=stealth,line width=1.2pt]
    (A)--(H1) node[left=3pt] {$h_1$};
  \draw[->,>=stealth,line width=1.2pt]
    (A)--(H2) node[right=3pt] {$h_2$};
  \draw[->,>=stealth,line width=1.2pt]
    (A)--(H3) node[below=4pt] {$h_3$};
\end{scope}

\fill (A) circle (1pt);
\fill (L1) circle (1pt);
\fill (L2) circle (1pt);
\fill (L12) circle (1pt);
\node[below left] at (A) {$A$};
\node[below=3pt] at (B) {$B$};
\node[right=3pt] at (C) {$C$};
\end{tikzpicture}
\caption{The reflection tiling of the plane.}
\label{fig:triangle-tiling}
\end{figure}

Translations are orientation-preserving, so they act with positive
sign on both the even and odd reflection extensions. Thus both
extensions are $\Lambda$-periodic and may be viewed as functions on
the flat torus 
$ \TT=\R^2/\Lambda.
$
See Figure \ref{fig:triangle-tiling}.

Neumann and Dirichlet eigenfunctions on $T$ extend to torus
eigenfunctions that are, respectively, even and odd under the
side reflections. Reflection extension naturally identifies $L^2(T)$ with the
corresponding reflection subspace of $L^2(\TT)$. We shall use this identification when
transferring operators from the torus to the triangle.

Set
\[
 h_1=\frac{t_2}{4}
 =\left(\frac{\sqrt3}{8},\frac38\right),\quad
 h_2=\frac{t_1-t_2}{4}
 =\left(\frac{\sqrt3}{8},-\frac38\right),\quad
 h_3=\frac{t_1}{4}
 =\left(\frac{\sqrt3}{4},0\right).
\]
The translation operators $T_dF(z)=F(z+d)$, $d\in\TT$, are
unitary on $L^2(\TT)$, commute with the torus Laplacian, and
satisfy $T_d^*=T_{-d}$. Each side reflection
$r(z)=Rz+a$ preserves the lattice through its linear part $R$, and
its pullback $r^*F=F\circ r$ satisfies
\begin{equation}\label{eq:triangle-translation-conjugation}
 r^*T_d(r^*)^{-1}=T_{R^{-1}d}.
\end{equation}

\subsection{The Neumann pair}\label{sec:neumann}

We prove Theorem \ref{thm:neumann}. Let $u$ be a finite sum of
Neumann eigenfunctions, and denote its even reflection extension by
$\widetilde u$.

Consider the subgroup of half-periods
\[
 K:=\{k\in\TT:2k=0\}
 =\{0,2h_1,2h_2,2h_3\}
 \cong(\Z/2\Z)^2.
\]
Let
\[
 \cP_K=\frac14\sum_{k\in K}T_k
\]
be the orthogonal projection onto the $K$-invariant functions on
$\TT$, and set
\begin{equation}\label{eq:neumann-torus-reflection}
 \cS=2\cP_K-I.
\end{equation}
Since the linear parts of the side reflections preserve $\Lambda$,
they also preserve $K$. By \eqref{eq:triangle-translation-conjugation}, $\cP_K$ and $\cS$
preserve the even reflection subspace and hence induce operators
on $L^2(T)$, denoted by the same symbols. Since $\cP_K$ and $\cS$ commute with the torus Laplacian and preserve
the even reflection subspace, their induced operators preserve
every Neumann eigenspace. 
Moreover, $\cP_K$ is an orthogonal projection, so $\cS$ is an
orthogonal involution.

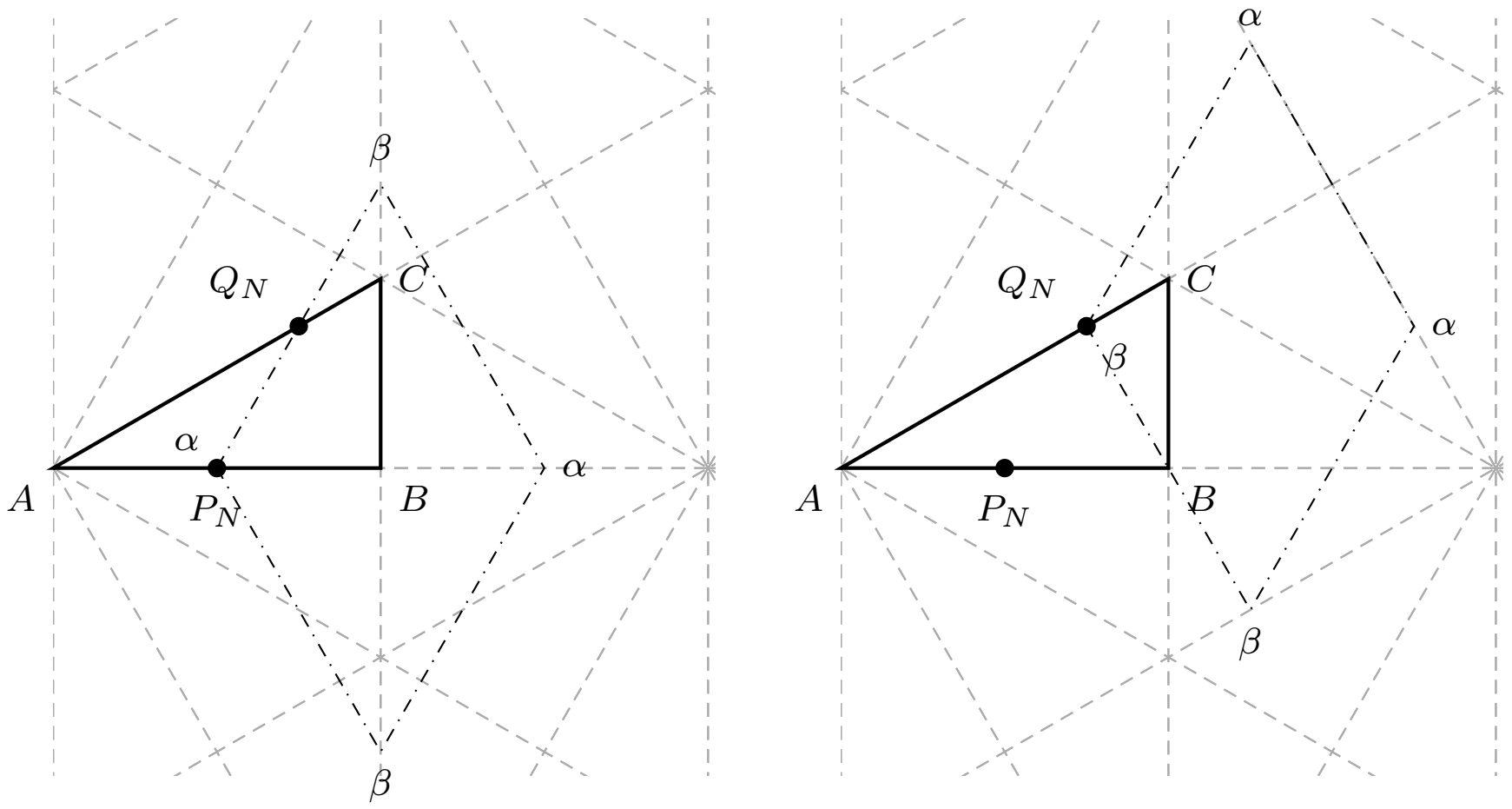
\begin{figure}[htbp]
\centering

\newcommand{\TriangleTiling}{
  \begin{scope}
    \clip (0,-13) rectangle (28,19);
    \foreach \tilei in {-1,0,1} {
      \foreach \tilej in {-1,0,1} {
        \begin{scope}[
          shift={({(16*\tilei+8*\tilej)*sqrt(3)},{24*\tilej})}]
          \foreach \tileangle in {0,60,120,180,240,300} {
            \begin{scope}[rotate=\tileangle]
              \draw[densely dashed,gray!65]
                (0,0)--({8*sqrt(3)},0)--
                ({8*sqrt(3)},8)--cycle;
              \draw[densely dashed,gray!65]
                (0,0)--({8*sqrt(3)},0)--
                ({8*sqrt(3)},-8)--cycle;
            \end{scope}
          }
        \end{scope}
      }
    }
  \end{scope}
  \draw[thick]
    (0,0)--({8*sqrt(3)},0)--({8*sqrt(3)},8)--cycle;
  \node[below left] at (0,0) {$A$};
  \node[below right] at ({8*sqrt(3)},0) {$B$};
  \node[right] at ({8*sqrt(3)},8) {$C$};
}

\begin{tikzpicture}[x=.18cm,y=.18cm,font=\scriptsize]

\begin{scope}
  \TriangleTiling

  \draw[dash dot]
    ({4*sqrt(3)},0)--({8*sqrt(3)},12)--
    ({12*sqrt(3)},0)--({8*sqrt(3)},-12)--cycle;

  \fill ({4*sqrt(3)},0) circle (2pt);
  \fill ({6*sqrt(3)},6) circle (2pt);
  \node[below=2pt] at ({4*sqrt(3)},0) {$P_N$};
  \node[above left=2pt] at ({6*sqrt(3)},6) {$Q_N$};

  \node[above left] at ({4*sqrt(3)},0) {$\alpha$};
  \node[above] at ({8*sqrt(3)},12) {$\beta$};
  \node[right] at ({12*sqrt(3)},0) {$\alpha$};
  \node[below] at ({8*sqrt(3)},-12) {$\beta$};

  \node at (14,-16) {The four translates of $P_N$};
\end{scope}

\begin{scope}[xshift=6cm]
  \TriangleTiling

  \draw[dash dot]
    ({6*sqrt(3)},6)--({10*sqrt(3)},18)--
    ({14*sqrt(3)},6)--({10*sqrt(3)},-6)--cycle;

  \fill ({4*sqrt(3)},0) circle (2pt);
  \fill ({6*sqrt(3)},6) circle (2pt);
  \node[below=2pt] at ({4*sqrt(3)},0) {$P_N$};
  \node[above left=2pt] at ({6*sqrt(3)},6) {$Q_N$};

  \node[below right] at ({6*sqrt(3)},6) {$\beta$};
  \node[above] at ({10*sqrt(3)},18) {$\alpha$};
  \node[right] at ({14*sqrt(3)},6) {$\alpha$};
  \node[below] at ({10*sqrt(3)},-6) {$\beta$};

  \node at (14,-16) {The four translates of $Q_N$};
\end{scope}

\end{tikzpicture}
\caption{The Neumann pair in the reflected triangle tiling.}
\label{fig:neumann-reflections}
\end{figure}

It remains to check what $\cS$ does at the marked points. Write
\[
 \alpha=u(P_N),\qquad \beta=u(Q_N).
\]
The four translates of $P_N$ and $Q_N$ are shown in
Figure \ref{fig:neumann-reflections}. One sees that
\[
\begin{array}{c|rrrr}
 k&0&2h_3&2h_1&2h_2\\ \hline
 \widetilde u(P_N+k)&\alpha&\alpha&\beta&\beta\\
 \widetilde u(Q_N+k)&\beta&\alpha&\alpha&\beta
\end{array}
\]
and therefore
\[
 (\cP_Ku)(P_N)=(\cP_Ku)(Q_N)=\frac{\alpha+\beta}{2}.
\]
It follows from \eqref{eq:neumann-torus-reflection} that
\begin{equation}\label{eq:neumann-exchange}
 (\cS u)(P_N)=u(Q_N),\qquad
 (\cS u)(Q_N)=u(P_N).
\end{equation}
Thus $\cS$ is a spectral symmetry exchanging $P_N$ and $Q_N$.
Theorem \ref{thm:spectral-symmetry} shows that the two points are
cospectral. Since $T$ is scalene, no isometry of $T$ exchanges them.
This proves Theorem \ref{thm:neumann}.

\subsection{The Dirichlet pair}\label{sec:dirichlet}
We give an alternative proof of the following result.

\begin{theorem}[\cite{DMS26}]
\label{thm:dirichlet}
The points $P_D,Q_D$ defined in \eqref{eq:dirichlet-points}
are cospectral for the Dirichlet Laplacian on $T$.
\end{theorem}

\begin{figure}[htbp]
\centering
\begin{tikzpicture}[x=.19cm,y=.19cm,font=\scriptsize]

\begin{scope}

\begin{scope}
  \clip (-1,-4.5) rectangle (21,12);
  \foreach \i in {-1,0,1} {
    \foreach \j in {-1,0,1} {
      \begin{scope}[
        shift={({(16*\i+8*\j)*sqrt(3)},{24*\j})}]
        \foreach \ang in {0,60,120,180,240,300} {
          \begin{scope}[rotate=\ang]
            \draw[densely dashed,gray!40,line width=.3pt]
              (0,0)--({8*sqrt(3)},0)--
              ({8*sqrt(3)},8)--cycle;
            \draw[densely dashed,gray!40,line width=.3pt]
              (0,0)--({8*sqrt(3)},0)--
              ({8*sqrt(3)},-8)--cycle;
          \end{scope}
        }
      \end{scope}
    }
  }
\end{scope}

\draw[thick]
  (0,0)--({8*sqrt(3)},0)--({8*sqrt(3)},8)--cycle;
\draw[densely dashed]
  (0,0)--({4*sqrt(3)},12)--({8*sqrt(3)},8);
\draw[densely dashed]
  (0,0)--({4*sqrt(3)},-4);

\draw[dash dot]
  ({9*sqrt(3)},3)--({7*sqrt(3)},9)--
  ({3*sqrt(3)},9)--({sqrt(3)},3)--
  ({3*sqrt(3)},-3)--({7*sqrt(3)},-3)--cycle;

\fill ({5*sqrt(3)},3) circle (2pt);
\fill ({7*sqrt(3)},3) circle (2pt);

\node[below left] at (0,0) {$A$};
\node[below right] at ({8*sqrt(3)},0) {$B$};
\node[right] at ({8*sqrt(3)},8) {$C$};

\node[below left] at ({5*sqrt(3)},3) {$Q_D$};
\node[above] at ({7*sqrt(3)},3) {$P_D$};

\node[left] at ({sqrt(3)},3) {$0$};
\node[above] at ({3*sqrt(3)},9) {$0$};
\node[above] at ({7*sqrt(3)},9) {$0$};
\node[below] at ({3*sqrt(3)},-3) {$0$};
\node[right] at ({9*sqrt(3)},3) {$-\alpha$};
\node[below] at ({7*sqrt(3)},-3) {$-\alpha$};

\node at ({5*sqrt(3)},-6)
  {Hexagon centered at $Q_D$};

\end{scope}

\begin{scope}[xshift=5.7cm]

\begin{scope}
  \clip (-1,-4.5) rectangle (21,12);
  \foreach \i in {-1,0,1} {
    \foreach \j in {-1,0,1} {
      \begin{scope}[
        shift={({(16*\i+8*\j)*sqrt(3)},{24*\j})}]
        \foreach \ang in {0,60,120,180,240,300} {
          \begin{scope}[rotate=\ang]
            \draw[densely dashed,gray!40,line width=.3pt]
              (0,0)--({8*sqrt(3)},0)--
              ({8*sqrt(3)},8)--cycle;
            \draw[densely dashed,gray!40,line width=.3pt]
              (0,0)--({8*sqrt(3)},0)--
              ({8*sqrt(3)},-8)--cycle;
          \end{scope}
        }
      \end{scope}
    }
  }
\end{scope}

\draw[thick]
  (0,0)--({8*sqrt(3)},0)--({8*sqrt(3)},8)--cycle;
\draw[densely dashed]
  ({8*sqrt(3)},8)--({10*sqrt(3)},10);

\draw[dash dot]
  ({11*sqrt(3)},3)--({9*sqrt(3)},9)--
  ({5*sqrt(3)},9)--({3*sqrt(3)},3)--
  ({5*sqrt(3)},-3)--({9*sqrt(3)},-3)--cycle;

\fill ({5*sqrt(3)},3) circle (2pt);
\fill ({7*sqrt(3)},3) circle (2pt);

\node[below left] at (0,0) {$A$};
\node[below right] at ({8*sqrt(3)},0) {$B$};
\node[right] at ({8*sqrt(3)},8) {$C$};

\node[below left] at ({5*sqrt(3)},3) {$Q_D$};
\node[above] at ({7*sqrt(3)},3) {$P_D$};

\node[left] at ({3*sqrt(3)},3) {$0$};
\node[above] at ({9*sqrt(3)},9) {$0$};
\node[above] at ({5*sqrt(3)},9) {$-\alpha$};
\node[right] at ({11*sqrt(3)},3) {$-\beta$};
\node[below] at ({5*sqrt(3)},-3) {$-\beta$};
\node[below] at ({9*sqrt(3)},-3) {$\alpha$};

\node at ({7*sqrt(3)},-6)
  {Hexagon centered at $P_D$};

\end{scope}

\end{tikzpicture}
\caption{The Dirichlet hexagons.}
\label{fig:dirichlet}
\end{figure}
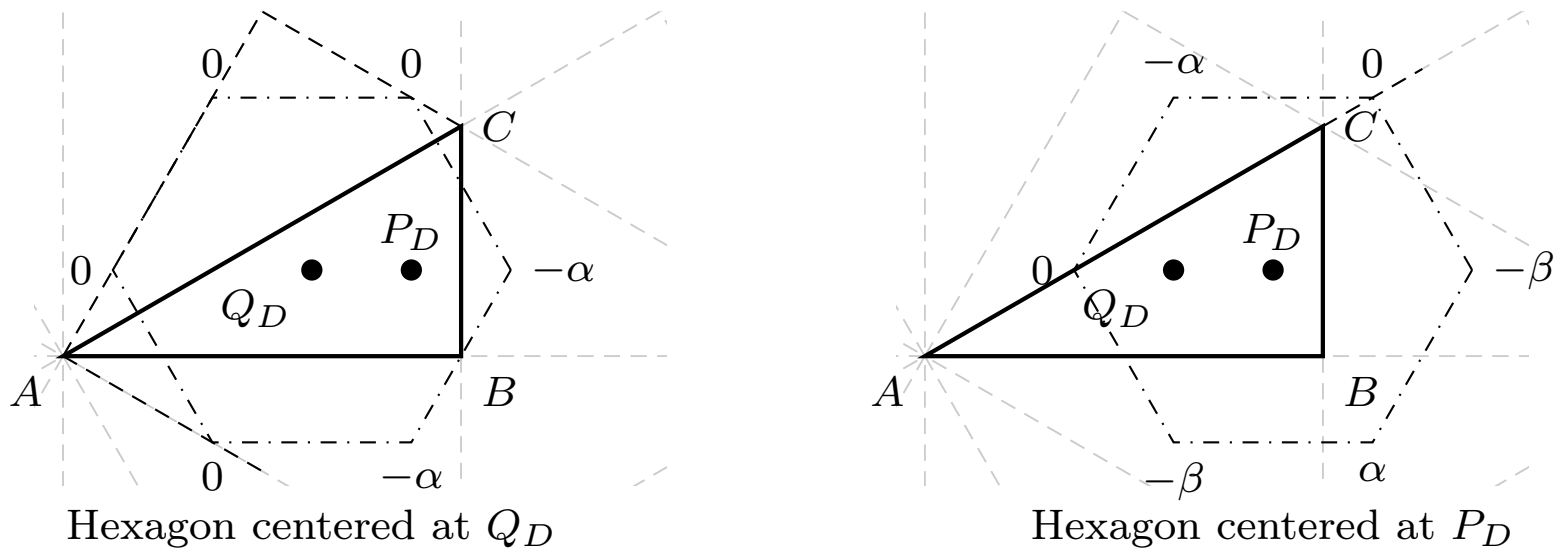

\begin{proof}
Let $u$ be a finite sum of Dirichlet eigenfunctions, and denote its
odd reflection extension by $\widetilde u$. Thus $\widetilde u$
changes sign under reflection in every mirror of the tiling and
vanishes on the mirror lines.

Let
\[
 \cD=\{\pm h_1,\pm h_2,\pm h_3\},
\]
whose elements are the vertices of a regular hexagon centered at
the origin. On $L^2(\TT)$, define
\[
 \cA_{\cD}=-\frac12\sum_{d\in\cD}T_d.
\]
The linear parts of the side reflections permute $\cD$, so
\eqref{eq:triangle-translation-conjugation} shows that $\cA_{\cD}$
preserves the odd reflection subspace. It therefore induces an
operator on $L^2(T)$, denoted by the same symbol. Since $\cA_{\cD}$
commutes with the torus Laplacian and preserves the odd reflection
subspace, its induced operator preserves every Dirichlet eigenspace.
Moreover, $T_d^*=T_{-d}$ and $\cD=-\cD$, so $\cA_{\cD}$ is bounded
and self-adjoint.

In terms of the odd extension $\widetilde u$, we have
\begin{equation}\label{eq:six-on-drum}
 (\cA_{\cD}u)(z)
 =-\frac12\sum_{d\in\cD}\widetilde u(z+d),
 \qquad z\in T.
\end{equation}

It remains to evaluate $\cA_{\cD}$ at the marked points. Write
$\alpha=u(P_D)$ and $\beta=u(Q_D)$. Folding the six vertices of
each hexagon back into $T$ gives
\[
\begin{array}{c|rrrrrr}
 d&h_3&h_1&-h_2&-h_3&-h_1&h_2\\ \hline
 \widetilde u(Q_D+d)&-\alpha&0&0&0&0&-\alpha\\
 \widetilde u(P_D+d)&-\beta&0&-\alpha&0&-\beta&\alpha
\end{array}
\]
as illustrated in Figure \ref{fig:dirichlet}. For example,
$Q_D+h_3$ is the reflection of $P_D$ across the vertical side of
$T$, so its value is $-\alpha$. The zero entries correspond to
points on mirror lines.

Multiplying the row sums by $-1/2$ gives
\begin{equation}\label{eq:dirichlet-auxiliary-exchange}
 (\cA_{\cD}u)(P_D)=u(Q_D),\qquad
 (\cA_{\cD}u)(Q_D)=u(P_D).
\end{equation}
Lemma \ref{lem:exchange} now proves the theorem.
\end{proof}

\subsection{Strong cospectrality and perfect state transfer}
\label{sec:triangle-strong-cospectrality}

Following the terminology for graphs in \cite{GS17}, two points
$x,y$ are \emph{strongly cospectral} if on each Laplace eigenspace
$E_\mu$ there is a sign $\varepsilon_\mu\in\{\pm1\}$ such that
\[
 u(y)=\varepsilon_\mu u(x)
 \qquad\text{for every }u\in E_\mu.
\]

\begin{proposition}\label{prop:triangle-strong-cospectrality}
Both triangle pairs are strongly cospectral and exhibit perfect
Schr\"odinger state transfer:
\begin{equation}\label{eq:triangle-state-transfer}
 e^{i\frac{9}{16\pi}\Delta_N}\delta_{P_N}=\delta_{Q_N},
 \qquad
 e^{i\frac{9}{32\pi}\Delta_D}\delta_{P_D}=i\,\delta_{Q_D}.
\end{equation}
Here $\Delta_N,\Delta_D$ denote the Neumann and Dirichlet
Laplacians, respectively. The identities are understood through
spectral expansions, with $\delta_x$ denoting point evaluation.
\end{proposition}
\begin{proof}
The eigenfunctions on $\TT$, indexed by $(m,n)\in\Z^2$,
have eigenvalues
\begin{equation}\label{eq:triangle-q}
 \frac{16\pi^2}{9}q,\qquad q=m^2-mn+n^2.
\end{equation}
Since $q$ is even exactly when $m,n$ are both even, the definition
of $\cP_K$ shows that $\cS=2\cP_K-I$ has multiplier $(-1)^q$.
The exchange identity \eqref{eq:neumann-exchange} gives
\begin{equation}\label{eq:strong-neumann}
 u(Q_N)=(-1)^q u(P_N).
\end{equation}

For odd $q$, the multiplier of $\cA_{\cD}$ is
\begin{equation}\label{eq:dirichlet-multiplier}
 -\left(
 \cos\frac{\pi m}{2}
 +\cos\frac{\pi n}{2}
 +\cos\frac{\pi(m-n)}{2}
 \right)=(-1)^{(q+1)/2}.
\end{equation}
Thus \eqref{eq:dirichlet-auxiliary-exchange} gives
\begin{equation}\label{eq:strong-dirichlet}
 u(Q_D)=(-1)^{(q+1)/2}u(P_D).
\end{equation}
For even $q$, the extension is $K$-invariant, so opposite vertices
of each hexagon have equal values. The Dirichlet folding table
then forces $u(P_D)=u(Q_D)=0$. This proves strong cospectrality.

At the stated times, the Schr\"odinger multipliers are
$e^{-i\pi q}=(-1)^q$ in the Neumann case and
$e^{-i\pi q/2}=i(-1)^{(q+1)/2}$ for odd $q$ in the Dirichlet case.
Together with the vanishing for even $q$, these identities give
\eqref{eq:triangle-state-transfer}.
\end{proof}

\section{A slit domain}\label{sec:slit}

Finally, we give a planar Neumann example with both marked points in the
interior. Recall
\[
 \Omega=\bigl((-1,1)\times(0,3)\bigr)\setminus
 \Bigl(\{0\}\times\Bigl[0,\frac12\Bigr]\Bigr),
\]
and the two points are
\[
 p=(0,1),\qquad q=(0,2).
\]

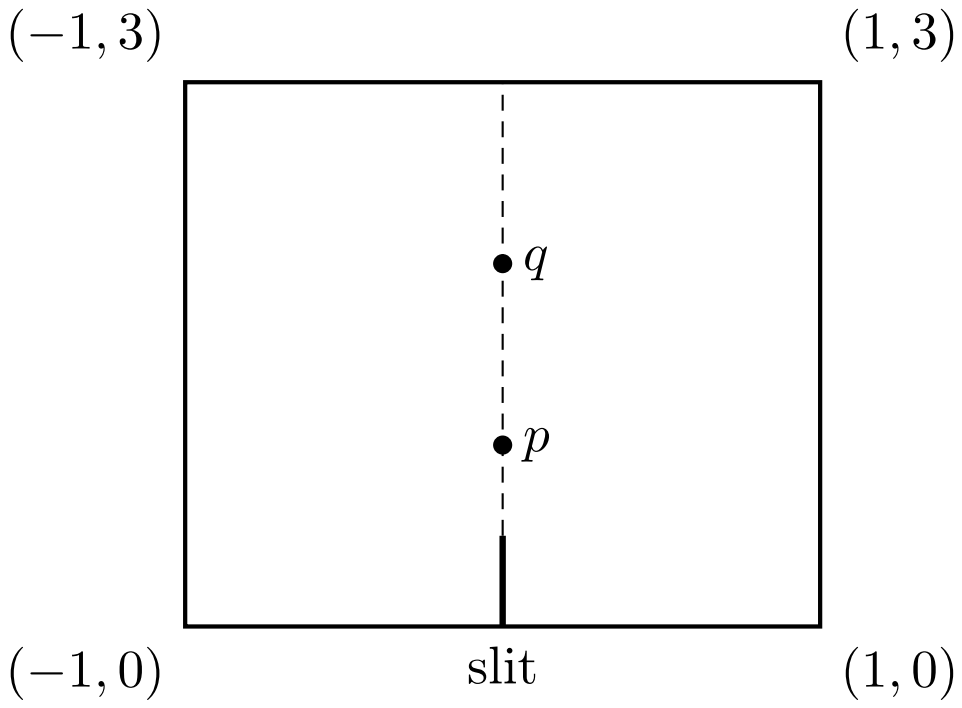
\begin{figure}[htbp]
\centering
\begin{tikzpicture}[x=2.1cm,y=1.2cm,font=\small]
  \draw[thick] (-1,0) rectangle (1,3);

  \draw[line width=1.2pt] (0,0)--(0,.5);

  \draw[densely dashed] (0,.5)--(0,3);

  \fill (0,1) circle (1.8pt);
  \fill (0,2) circle (1.8pt);

  \node[right] at (0,1) {$p$};
  \node[right] at (0,2) {$q$};
  \node[below] at (0,0) {slit};

  \node[below left] at (-1,0) {$(-1,0)$};
  \node[below right] at (1,0) {$(1,0)$};
  \node[above left] at (-1,3) {$(-1,3)$};
  \node[above right] at (1,3) {$(1,3)$};
\end{tikzpicture}
\caption{The slit domain $\Omega$. }
\label{fig:slit}
\end{figure}

\begin{proof}[Proof of Theorem \ref{thm:slit}]
Let $\sigma(x,y)=(-x,y)$. Since $\sigma$ is an isometry of $\Omega$,
the Neumann Laplacian preserves the decomposition of $L^2(\Omega)$
into the even and odd subspaces
\[
 L^2(\Omega)=\mathcal H_+\oplus\mathcal H_-.
\]

The odd part is invisible at the marked points. Indeed, every Neumann eigenfunction 
$u\in\mathcal H_-$ vanishes on the fixed axis $x=0$, and hence
$ u(p)=u(q)=0.
$

For the even part, the slit is invisible. An even Neumann eigenfunction
has the same trace on the two faces of the slit, so the two sides glue
together and the function extends across the slit to a Neumann
eigenfunction on the unslit rectangle
$ (-1,1)\times(0,3).
$
On this rectangle, reflection across the horizontal midline,
$\rho(x,y)=(x,3-y)$, preserves the Neumann Laplacian and exchanges
$p$ and $q$. Since $\rho$ also preserves evenness in $x$, extending
across the slit, reflecting by $\rho$, and restricting back to
$\Omega$ defines an orthogonal involution
$J:\mathcal H_+\to\mathcal H_+$ preserving every Neumann eigenspace,
with
\[
 (Ju)(p)=u(q),\qquad (Ju)(q)=u(p).
\]

Now define $S$ on $L^2(\Omega)=\mathcal H_+\oplus\mathcal H_-$ by
\[
 S(u_++u_-)=Ju_++u_-.
\]
Then $S$ is an orthogonal involution preserving every Neumann
eigenspace. If $u=u_++u_-$ is a finite sum of Neumann eigenfunctions,
the odd part vanishes at both marked points, so
\[
 (Su)(p)=(Ju_+)(p)=u_+(q)=u(q),
\]
and similarly $(Su)(q)=u(p)$. Thus $S$ is a spectral symmetry
exchanging $p$ and $q$. Theorem \ref{thm:spectral-symmetry} therefore
shows that $p$ and $q$ are cospectral.

Finally, any isometry of $\Omega$ must preserve the outer rectangle
and the slit. The slit rules out the horizontal reflection and the
half-turn, leaving only the identity and reflection in the vertical
axis. Both fix $p$ and $q$ individually, so no isometry of $\Omega$
exchanges them.
\end{proof}

This example fits especially directly into the spectral-symmetry
framework, but by a different mechanism from the averaging and
unfolding constructions. The vertical reflection does not exchange the
two points. Instead, it splits the spectral problem into even and odd
parts. The odd part is invisible at the marked points, while the even
part reduces to a half-rectangle with a horizontal reflection exchanging
them.
\bibliographystyle{amsalpha}
\bibliography{reference}
\end{document}